\documentclass[12pt]{article}
\usepackage{amsmath,amssymb}
\newtheorem{thm}{Theorem}[section]

\newtheorem{cor}{Corollary}[section]
\newtheorem{prop}{Proposition}[section]

\newcommand{\x}{\times}

\newcommand{\bs}{\bigskip}

\newcommand{\ds}{\displaystyle}
\renewcommand{\d}{\delta}

\newcommand{\M}{\mathcal M}
\newcommand{\mb}{\mbox}

\renewcommand{\O}{\mathcal O}

\newcommand{{\z}}{\Bbb Z}
\renewcommand{\P}{\mathbb P}
\newcommand{\ci}{Clifford index }
\newcommand{\gc}{geometric}
\newcommand{\go}{gonality }
\newcommand{\emb}{embedding }
\newcommand{\crg}{curves of genus }
\newcommand{\hy}{hyperelliptic }
\newcommand{\gon}{\mathop{\mathrm{gon}}\nolimits}
\newcommand{\de}{\mathop{\mathrm{deg}}\nolimits}

\begin{document}
\setlength{\baselineskip}{18pt}

\begin{center}
{\Large\bf Clifford Index of ACM Curves in ${\mathbb
P}^3$}\bs\bs\\
{\sc robin hartshorne}\bs\\
Department of Mathematics\\University of California\\
Berkeley, California 94720--3840
\end{center}

\bs\bs
\begin{quote}
In this paper we review the notions of gonality and
Clifford index of an abstract curve. For a curve embedded
in a projective space, we investigate the connection
between the \ci of the curve and the \gc al properties of
its \emb. In particular if $C$ is a curve of degree $d$
in ${\P}^3$, and if $L$ is a multisecant of maximum order
$k$, then the pencil of planes through $L$ cuts out a
$g^1_{d-k}$ on $C$. If the gonality of $C$ is equal to
$d-k$ we say the gonality of $C$ can be computed by
multisecants. We discuss the question whether the \go of
every smooth ACM curve in ${\P}^3$  can be computed by
multisecants, and we show the answer is yes in some
special cases.
\end{quote}

\section{Gonality and \ci of a curve}

Let $C$ be a nonsingular projective curve over an
algebraically closed field $k$. A linear system of
degree $d$ and (projective) dimension $r$ will be denoted by
$g^r_d$. The least integer $d$ for which there exists a complete
linear system $g^1_d$ without base points is called the
{\it \go} of $C$. Thus a curve is rational if and only if its \go
is 1. Curves of genus 1 and 2 have \go 2. For curves of genus
$g\geq 2$, the curve is hyperelliptic if and only if the \go is
2. It is well known that for \crg $g\geq 3$ the \go $d$ lies
between 2 and $\ds{\left[\frac{g+3}{2}\right]}$\,;
there exist curves having each possible \go in this range; and a
\crg $g$ of general moduli has \go
$\ds{\left[\frac{g+3}{2}\right]}$\,.
See \cite{ACGH}  for references to proofs of these
results.

Thus the \go of a curve provides a stratification of the variety
of moduli ${\M}_g$ of \crg $g$, with the hyperelliptic curves at
one end, and the curves of general moduli at the other end.

To illustrate this principle, let us describe some different types
of \crg $g$ for small values of $g$.

For $g=0$ there is just one curve, ${\P}^1$, having \go 1 and a
unique $g^1_1$.

For $g=1$ there is a one-parameter family of non-isomorphic
curves. They all have \go 2, and each one has infinitely many
$g^1_2$'s.

For $g=2$, the curves are \hy\!\!, each having a unique $g^1_2$.

For $g=3$, there are \hy curves, with a unique $g^1_2$, and there
are non-\hy curves, each having infinitely many   $g^1_3$'s.
These are called {\it trigonal} (meaning a curve with \go 3).
The canonical \emb of a trigonal curve of genus 3 is a nonsingular
plane quartic curve. The $g^1_3$'s on the curve are cut out by the
pencils of lines through a point on the curve.

For $g=4$, there are again two types, \hy and trigonal. The
canonical \emb of a trigonal \crg 4 is a complete intersection of
an irreducible quadric surface $Q$ and a cubic surface $F_3$ in
${\P}^3$. If $Q$ is nonsingular (the general case), the the curve
$C$ has exactly two $g^1_3$'s cut out by the two families of lines
on $Q$.  If $Q$ is a cone, then $C$ has a unique $g^1_3$.

For $g=5$ there are three types of curves: the \hy curves, the
trigonal curves, and the general curves. The canonical \emb of a
non-\hy curve of 5 is a curve of degree 8 in ${\P}^4$.
If the curve is trigonal, then it lies on a rational ruled cubic
surface in ${\P}^4$, and the unique $g^1_3$ is cut out by the
rulings of that surface. In the general case, the curve is a
complete intersection of three quadric hypersurfaces  in
${\P}^4$, and has infinitely many $g^1_4$'s.

For $g=6$ the situation becomes more complicated and more
interesting. We can distinguish (at least) five different types of
curves:
\begin{description}
\item{a)}  The \hy curves, having a unique $g^1_2$.
\item{b)} The trigonal curves, having a unique $g^1_3$.
\item{c)} Plane quartic curves, having a unique $g^2_5$.
These curves have infinitely many $g^1_4$'s, cut out by the pencils
of lines through a point of the curve.
\item{d)} Double cover of an elliptic curve,
having infinitely many $g^1_4$'s, but no $g^2_5$.
\item{e)} Curves having only finitely many $g^1_4$'s.
The general curve has exactly five $g^1_4$'s; some others may have
fewer.
\end{description}

For references, see e.g.~\cite[Ch.V]{ACGH}.

   From these few examples it is already clear that the \go does not
tell the whole story in distinguishing different types of \crg
$g$. More generally, we should take into account all possible
special linear systems $g^r_d$ that might exist on the curve.
Here {\it special} means $r>d-g$, or equivalently
$H^1(C,{\O}_C(D))\neq 0$, where $D$ is a divisor in the linear
system.

In this connection we consider the {\it Brill-Noether number}
\[
\rho = g-(r+1)(g-d+r) \ .
\]
Then one knows, for given $g,d,r$ that if $\rho\geq 0$,
every \crg $g$ has a $g^r_d$, and if $\rho<0$, then a general \crg
$g$ does not have a $g^r_d$ \cite[Ch.V]{ACGH}.
So to distinguish among \crg $g$, we will be interested in the
existence of linear systems $g^r_d$ for which $\rho<0$.

Now we can introduce the {\it \ci} of a curve $C$ of genus $g\geq
4$ (first defined by H.~Martens \cite{HM}). For a particular linear
system $g^r_d$ we define its \ci to be $d-2r$.
Then the {\it \ci of the curve} is the minimum of Clifford indices
of certain special linear series, namely Cliff $C$ is the minimum
of $d-2r$, taken over all linear systems $g^r_d$ with $r\geq 1$
and $1\leq d\leq g\!-\!1$. Equivalently, one can take the minimum
of $d-2r$ over all $g^r_d$ containing a divisor $D$ with
$h^0({\O}_C(D))\geq 2$ and $h^1(({\O}_C(D))\geq 2$.
(The equivalence of these two criteria is easy using the
Riemann-Roch theorem, and replacing $D$ by $K-D$ if $d>g-1$.)

Recall that Clifford's theorem tells us that $r\leq \frac 12 d$
for a special linear system $g^r_d$ on a curve, with equality if
and only if the corresponding divisor $D$ is 0 or the canonical
divisor $K$, or the curve is \hy\!\!. Since the inequalities in
the definition rule out the possibility $D=0$ or $K$, we see that
the Clifford index is always $\geq 0$, with equality if and only
if the curve is \hy\!\!. On the other hand, the Brill-Noether
theory tells us that Cliff
$C\leq\ds{\left[\frac{d-1}{2}\right]}$\,, and is equal to this
value for a curve of general moduli.

In most cases, the \ci can be computed by a pencil, that is,
there exists a $g^1_d$ with Cliff $C=d-2$. In this case
Cliff $C={\mb{gon}}\,C-2$, where ${\mb{gon}}\,C$ denotes the gonality.
This suggests the definition of the {\it Clifford dimension} of the
curve, which is the least $r$ for which there exists a $g^r_d$ with
Cliff
$C=d-2r$, i.e., the $g^r_d$ computes the \ci\!\!. Then $r=1$ is
the normal situation. Curves with Clifford dimension $>1$ are
rare. The first example is the curve of genus 6 with a $g^2_5$
mentioned above. In this case Cliff $C=1$ while ${\mb{gon}}\,C=4$, and so the
Clifford dimension  of the $g^1_4$ is 2.

The nonsingular plane curves of degree $\geq 5$ all have
Clifford dimension 2, and these are the only such. Any curve of
Clifford dimension 3 must be a complete intersection of two cubic
surfaces in ${\P}^3$, having degree 9 and genus 10 \cite{GM1}.
There exist curves of every possible Clifford dimension $r\geq 1$,
and for $r\geq 3$ conjecturally only one possible degree-genus
pair in ${\P}^r$ \cite{ELMS}.

\section{Curves in projective space}

We now consider a nonsingular curve $C$ embedded in a projective
space ${\P}^n_{\!k}$, and we ask, how are the \go and the \ci of
$C$ related to the geometry of the embedding? The prototype for
this kind of question is the following well-known theorem about
plane curves.

\begin{thm} Let $C\subseteq {\P}^2$ be a nonsingular plane curve
of degree $d\geq 2$. Then

{\rm{(a)}} There is no $g^1_{d-2}$ on $C$, but there are
$g^1_{d-1}$'s on $C$, so the \go is $d-1$.

{\rm{(b)}} Every $g^1_{d-1}$ on $C$ is cut out by the pencil of
lines in ${\P}^2$ through some fixed point of~$C$.
\end{thm}

{\bf Proof.} This result was known to M.Noether and has received a
number of modern proofs \cite{C},\cite{H},...  We will give an
elementary proof to illustrate the ideas involved.

(a) Suppose $D$ is a divisor of degree $e\leq d-2$ on $C$ with
$h^0({\O}(D))\geq 2$. Since $D$ moves in a pencil, we may assume
that $D$ consists of $e$ distinct points $P_1,\dots ,P_e$.
By the Riemann-Roch theorem $h^0(K-D)\geq g-e+1$. This means that
the $e$ points  $P_1,\dots ,P_e$ do not impose independent
conditions on the canonical divisors $K$ containing them.
Now the canonical divisor $K$ on $C$ is cut out by curves of
degree $d-3$ in ${\P}^2$. Any $d-2$ distinct points impose independent
conditions on these curves, a contradiction. So no such $D$ exists.
On the other hand, for any $P\in C$, the lines through $P$ cut out
a $g^1_{d-1}$, so these do exist.

(b) Now let $D$ be any divisor of degree $d-1$ with
$h^0({\O}(D))\geq 2$. The argument above shows that the $d-1$
points of $D$ impose dependent conditions on plane curves of
degree $d-3$, and this can only happen if these points lie on a
line in ${\P}^2$.  This line $L$ will meet $C$ in one further
point $P$, and then it is clear that the $g^1_{d-1}$ is equal to
the one cut out by lines through $P$.

\bs This result has been generalized to irreducible plane curves
$C$ of degree $D$ with $\d$ nodes and cusps, when $\d$ is not too
large in relation to $d$ \cite{Co},\cite{CK}. In those cases, the
desired result would be that the \go of the normalization
$\tilde{C}$ is $d-2$, and is given by the linear systems cut out
by lines through one of the double points. We cannot expect such a
result to hold for arbitrary plane curves with nodes, however, as
the following example shows.

\bs{\bf Example 2.2.} Let $C$ be a smooth curve in ${\P}^3$, of
degree 6 and genus 3, not lying on a quadric surface. Such a curve
arises of type $(4;1^6)$ on a nonsingular cubic surface $X$ in
    ${\P}^3$, for example. This is the proper transform of a plane
curve of degree 4 passing through the 6 points blown up to get
$X$. It has \go 3 by $(2.1)$. On the other hand, its general
projection to  ${\P}^2$ is a plane curve of degree 6 with 7 nodes.
The pencil of lines through one of the nodes cuts out a $g^1_4$ on
the normalization, which does not give us the correct \go\!\!.

\bs
Passing now to curves in higher dimensional projective spaces, let
us first consider a nonsingular curve $C$ of degree $d$ in
    ${\P}^3$.  Let $L$ be a multisecant of maximum order $k$
(that is, the scheme-theoretic intersection of $C$ and $L$ has
length $k$). Then the pencil of planes through $L$ cuts out a
$g^1_{d-k}$ on $C$. If the \go of $C$ is equal to $d-k$, we say the
\go of $C$ can be {\it computed by multisecants.} We can also ask
the stronger question, whether every $g^1_{d-k}$ on $C$ arises in
this way.

For a curve $C$ of degree $d$ in ${\P}^n$, with $n\geq 4$,
the corresponding situation would be to look for a multisecant
linear space $L$ of codimension 2 in ${\P}^n$, meeting $C$ in a
scheme of length $k$. The hyperplanes through $L$ will cut out a
$g^1_{d-k}$ on $C$, and if this gives the \go\!\!, we say again
that the \go can be computed by  multisecants.

\bs{\bf Example 2.3.} Let $C\hookrightarrow{\P}^n$ be the
canonical \emb of a non\hy curve of genus $g\geq 3$ in ${\P}^n$,
with $n=g-1$. Let $g^1_e$ be a special complete linear system
without base points, and let $D$ be any divisor of the $g^1_e$.
Then $h^0({\O}(D))=2$. Let $F=K-D$, where $K$ is the canonical
divisor. Then  $h^0({\O}(K-F))=2$. Since $K$ is cut out by
hyperplanes in ${\P}^n$, this means that the divisor $F$ is
contained in two distinct hyperplanes. Let them meet in the linear
space $L$ of codimension 2. Then the pencil of hyperplanes through
$L$, after removing the fixed points $F$, cuts out the original
linear system $g^1_e$ on $C$. In particular, the \go can be
computed by multisecants.

If $g=3$, we have a plane curve of degree 4, and recover the
result of $(2.1)$ in this case.

If $g=4$, the curve $C$ is the complete intersection of a quadric
surface $Q$ with a cubic surface $F_3$. The pencil of planes
through a line of $Q$ cuts out the other family of lines on $Q$,
if $Q$ is nonsingular, or the only family of lines, if $Q$ is a
cone. Thus the $g^1_3$'s are computed by multisecants.

If $g=5$, there are two cases. When $C$ is trigonal, it lies on a
rational ruled cubic surface $S$ in ${\P}^4$. This surface
contains conics meeting the curve in 5 points. The plane of the
conic is therefore a 5-secant plane, and the hyperplane sections
of $S$ containing this conic cut out the rulings of $S$, which in
turn  cut out the unique $g^1_3$ on $C$.

If $C$ is not trigonal, then our result tells us that the curve
$C$ has 4-secant planes, so that the pencil of hypersurface
through them cut out the $g^1_4$'s on $C$.

\bs{\bf Example 2.4.} Basili \cite[$4.2$]{Bas} shows that if $C$ is a
smooth complete intersection curve in ${\P}^3$, not contained in a
plane, then the \go can be computed by multisecant lines.
Furthermore every $g^1_e$ with $e={\mb{gon}}\,C$ arises in this manner.

One can also ask what are the possible orders of multisecant
lines, and hence what are the possible gonalities of these complete
intersection curves. Let $C$ be the complete
intersection $F_a,F_b$ of surfaces of degrees $2\leq a\leq b$.
Nollet \cite{N} has shown that the maximum order $k$ of a
multisecant is either $\leq a$ or $=b$. If $a=2$, then $C$ is a
curve of bidegree $(b,b)$ on the quadric surface, so $k=b$ and
${\mb{gon}}\,C=d-b$. If
$a=3$, again
$k=b$ since there are lines on the cubic surface, so ${\mb{gon}}\,C=d-b$.
If $a\geq 4$, then Ellia and Franco \cite{EF} have shown that
every value of $k$ satisfying $4\leq k\leq a$ or $k=b$ can occur.
In particular, the general complete intersection curve with
$a,b\geq 4$ has at most 4-secants, and \go $d-4$.

\bs{\bf Example 2.5.} If $C$ is a smooth curve of bidegree $(a,b)$
on a nonsingular quadric surface, with $a\leq b$, then the maximum
order of a multisecant is $b$. In this case G.\ Martens \cite{GM1} and
Ballico \cite{Bal} have shown that the \go of $C$ is $a$, and thus is
computed by multisecants.

\bs{\bf Example 2.6.} The case of complete intersection curves in
${\P}^3$ $(2.4)$ has been generalized by Ellia and Franco \cite{EF} to
curves arising as the zero locus of a section of a rank 2 vector
bundle ${\mathcal E}$ on ${\P}^3$, twisted sufficiently:
$s\in H^0({\mathcal E}(t))$ with $t>>0$. This includes ``most"
subcanonical curves in ${\P}^3$, but it is hard to get specific
results for small degree curves. They show the \go of these curves
can also be computed by multisecants.

\bs{\bf Example 2.7.} Farkas \cite[$2.4$]{F} using the method of
Mori, shows the existence of curves $C$ on a nonsingular quartic
surface $X$ in ${\P}^3$, for certain values of degree $d$ and
genus $g$ satisfying complicated conditions, for which the \go is
$d-4$ and can be computed by 4-secants to the curve.

\bs{\bf Example 2.8.} Eisenbud et al \cite{ELMS}, also studying curves
on $K3$ surfaces, show the existence, for every $r\geq 3$ of a
smooth curve in ${\P}^r$ having $d=4r-3$, $g=4r-2$, whose \go is
computed by multisecants, and having Clifford dimension $r$.
The \go is $2r$ and the \ci is $2r-3$.

\bs
Lest the reader begin to think that all these examples are evidence
for supposing that the \go of any space curve is computed by
multisecants, let us give a few counterexamples.

\bs{\bf Example 2.9.} We consider rational curves $C$ of degree
$d$ in ${\P}^3$. There exist such curves having a
$(d\!-\!1)$-secant, for example, curves of bidegree $(1,d-1)$ on a
quadric surface. In this case the \go is computed by multisecants.
However, for $d\geq 5$, Ellia and Franco \cite{EF} have shown that
there exist smooth rational curves whose maximum order of a
multisecant $k$ can take any value $4\leq k\leq d-1$. In
particular, the general such curve has only 4-secants, and if
$d\geq 6$ these do not give the correct gonality.

\bs{\bf Example 2.10.} For an example of curves of higher genus,
take a plane curve of degree $a$, blow up 6 points not on the
curve, and let $C$ be the image curve in the nonsingular cubic
surface $X$ in ${\P}^3$. Then $C$ has degree $d=3a$. It has
multisecants of order $a$ and $2a$ on $X$. The pencil of planes
through one of the latter cut out a $g^1_a$ on $C$. But the \go is
$a-1$, so the \go is not computed by
multisecants.

\bs
   From all this evidence, it seems reasonable to pose the following
question.
\begin{quote} {\bf Question 2.11} (Peskine)
If $C$ is a smooth ACM (i.e. projectively normal)
curve in ${\P}^3$, is its \go  computable by
multisecants?\end{quote}

We will discuss this question in the following sections.

\section{Behavior of \go in a family}

We consider a smooth surface $X$, together with a proper morphism $f:
X \rightarrow T$,
where $T$ is a nonsingular curve.  Then the fibers of $f$ from a flat
family of curves
on $X$.  We assume that the general fiber $f^{-1}(t) = C_t$ is
irreducible and nonsingular for $t \in T$, $t \ne 0$, and that the
special fiber
$f^{-1}(0) = C_1 \cup
C_2$, a union of two smooth irreducible curves $C_1$, $C_2$, meeting
transversally at
$s$ distinct points.

\begin{thm} In a flat family as above, whose general curve $C_t$ is
smooth, and whose
special curve $C_0 = C_1 \cup C_2$ is a union of two smooth curves
meeting transversally
at $s$ points, we have
\[
\gon(C_t) \ge \min\{s,\gon(C_1) + \gon(C_2)\}
\]
for all sufficiently general $t \in T$.
\end{thm}

{\bf Proof.}  Suppose that the general curves $C_t$ in the family all
have a $g^1_d$ with
$d \le s$.  Then we can find an open set $T' \le T$ and an invertible
sheaf ${\mathcal
L}'$ on $X' = f^{-1}(T')$ inducing a $g^1_d$ on each fiber.  This
invertible sheaf
extends to an invertible sheaf ${\mathcal L}$ on $X$, but the
extension is not unique,
because we can replace ${\mathcal L}$ by ${\mathcal L}(mC_1)$ for any
$m \in {\mathbb
Z}$ and still get the same ${\mathcal L}'$ on restricting to $X'$.

Let us compute some intersection numbers.  Since $C_1 + C_2 \sim
C_t$, we have $C_1
\cdot (C_1 + C_2) = 0$, so $C_1^2 = -C_1\cdot C_2 = -s$.  Similarly
$C_2^2  = -s$.  Let us
denote by ${\mathcal L}_0$ the restriction of ${\mathcal L}$ to
$C_0$, and by ${\mathcal
L}_1,{\mathcal L}_2$, the restrictions to $C_1,C_2$.  Let $a = \de {\mathcal
L}_1$ and $b = \de {\mathcal L}_2$.  Then $a+b = \de {\mathcal L} = d$.
If we replace ${\mathcal L}$ by ${\mathcal L}(mC_1)$ then $a$ becomes
$a - ms$ while $b$
becomes $b+ms$.  Thus by choosing $m$ appropriately, we may assume
that $1 \le a\le s$
and consequently $1 - s \le b \le s-1$.

Since ${\mathcal L}$ cuts out a $g^1_d$ on the general curve $C_t$, we have
$h^0({\mathcal L}_t) \ge 2$ for general $t \in T$.  Hence by semicontinuity,
$h^0({\mathcal L}_0) \ge 2$.  We consider the exact sequence
\[
0 \rightarrow {\mathcal L}_0 \rightarrow {\mathcal L}_1 \oplus {\mathcal L}_2
\rightarrow {\mathcal O}_D \rightarrow 0
\]
where $D$ is the set of $s$ points $C_1 \cap C_2$.  On cohomology this gives
\[
0 \rightarrow H^0({\mathcal L}_0) \rightarrow H^0({\mathcal L}_1)
\oplus H^0({\mathcal
L}_2) \rightarrow H^0({\mathcal O}_D)
\]
we consider three cases, depending on the dimensions of $h^0({\mathcal L}_1)$,
$h^0({\mathcal L}_2)$.

\bs{\sc Case 1.}  If one of $h^0({\mathcal L}_i)$, $i = 1,2$ is zero,
then the other
must be $2$.  So we have a $g^1_d$ on one of the curves, say $C_1$,
and since the section
of ${\mathcal L}_0$ giving this $g^1_d$ is identically zero on $C_2$,
the divisor of the
$g^1_d$ must contain $D$ as a fixed component.  But then $d > s$,
contrary to our
assumptions.

\bs{\sc Case 2.}  If one of the $h^0({\mathcal L}_i)$ is $1$, say
$h^0({\mathcal L}_2) =
1$, then we can find a section of ${\mathcal L}_0$ inducing $0$ on
$C_2$ and a nonzero
section of ${\mathcal L}_1$ on $C_1$.  In this case $a \ge s$, and $b
\ge 0$ since
$h^0({\mathcal L}_2) \ne 0$, so $d \ge s$ and we have $\gon(C_t) \ge s$.

\bs{\sc Case 3.}  If both of $h^0({\mathcal L}_i)$ are $\ge 2$, then
we have a $g^1_a$ as
$C_1$ and a $g^1_b$ on $C_2$, so $d = a+b \ge \gon(C_1) + \gon(C_2)$, as
required.

\setcounter{cor}{1}
\begin{cor} In the statement of the theorem, if $\gon(C_t) =
\gon(C_1) + \gon(C_2) < s$,
then there exist morphisms $\pi_1: C_1 \rightarrow {\mathbb P}^1$ and
$\pi_2: C_2
\rightarrow {\mathbb P}^1$ of degrees equal to the gonality, such
that $\pi_1$ and
$\pi_2$ agree on the $s$ points $C_1 \cap C_2$.
\end{cor}

{\bf Proof.}  Indeed, if $\gon(C_t) < s$, then we must be in Case~3
of the proof above,
and the $g^1_a$ on $C_1$ and $g^1_b$ on $C_2$ are induced by
$H^0({\mathcal L}_0)$, so
must agree on $D$.

\bs{\it Note:}  A special case of this kind of argument appears in a
paper of Ballico
\cite{Bal}.

\bs{\bf Example 3.3.}  We can use the theorem to give another proof
of a weak form of
$(2.1)$, namely, a general curve $C$ of degree $d \ge 2$ in ${\mathbb
P}^2$ has gonality
$d-1$.  The pencil of lines through a point on $C$ cuts out a
$g^1_{d-1}$, so we always
have $\gon(C) \le d-1$.  To prove the reverse inequality, we use
induction on $d$.

For $d=2$, the conic is isomorphic to ${\mathbb P}^1$, so has gonality $1$.

For $d \ge 3$, consider a family of smooth curves $C_t$ of degree $d$
degenerating to
the union of a general smooth curve $C_1$ of degree $d-1$ and a
transversal line $C_2 =
L$.  Then $C_1 \cap L$ is $d-1$ points, and $\gon(C_1) = d-2$, $\gon
C_2 = 1$, so by the
theorem we find $\gon C \ge (d-2)+1 = d-1$.  Note that $s = \gon(C_1)
+ \gon(C_2)$ in
this proof, so we do not get any additional information from the Corollary.

\section{Curves on quadric and cubic surfaces}

\begin{prop} Let $Q$ be a nonsingular quadric surface in ${\mathbb
P}^3$, and let $C$ be
a general smooth curve of bidegree $(a,b)$ with $0 < a \le b$.  Then
$\gon(C) = a$.
\end{prop}

{\bf Proof.}  By projection onto one of the factors of $Q \cong
{\mathbb P}^1 \x
{\mathbb P}^1$ we know that $\gon(C) \le a$.  For $a=1$, the curve is
rational, so has
gonality $1$.

For $a \ge 2$, we let $C$ move in a family specializing to the union
of a general curve
$C_1$ of bidegree $(a-1,b-1)$ and a conic $C_2$ of bidegree $(1,1)$.
Then $s = C_1
\cdot C_2 = a+b-2 \ge a$.  Also by induction $\gon(C_1) = a-1$ and
$\gon(C_2) = 1$.  So
by $(3.1)$, $\gon(C) \ge a$, as required.

We recover a slightly weaker version of $(2.5)$, since our method of
proof applies only
to the general curve in a family.

\bs For curves on a cubic surface, we have seen $(2.10)$ that not
every smooth curve on
a smooth cubic surface has its gonality determined by multisecants.
However, we can
obtain a result for sufficiently general ACM curves on a cubic surface.

\begin{prop} Let $C$ be a smooth ACM curve on a nonsingular cubic
surface $X$ in
${\mathbb P}^3$.  If $C$ is sufficiently general in its linear system
on $X$, then $C$
has a multisecant $L$ such that the pencil of planes through $L$ cuts
out a pencil on
$C$ computing the gonality of $C$.
\end{prop}

{\bf Proof.}  First we must identify the smooth ACM curves on $X$.  Using the
postulation character $\gamma$ of \cite[$2.11$, p.~34]{MDP} we see
that if
$C$ is not contained in any quadric surface, then its
$\gamma$-character has $s_0 = 3$,
and is positive, and connected.  This means it must have one of the
following four types,
where
$a \ge 0$:
\begin{description}
\item{a)} \ $-1\ -1\ -1\ 0^a\ 3$
\item{b)} \ $-1\ -1\ -1\ 0^a\ 2 \ 1$
\item{c)}  \ $-1\ -1\ -1\ 0^a\ 1 \ 2$
\item{d)} \ $-1\ -1\ -1\ 0^a\ 1\ 1\ 1$.
\end{description}
These can all be obtained by ascending biliaison on $X$ from one of
the following curves
on a quadric surface:\\{}\\
\begin{tabular}{lll}
a) & $\gamma = -1\ -1\ 2$  &\quad $d = 3$, $g = 0$\\{}\\
b) & $\gamma = -1\ -1\ 1\ 1$ &\quad $d=4$, $g=1$\\{}\\
c) & $\gamma = -1-1\ 0 \ 2$ &\quad $d=5$, $g=2$\\{}\\
d) & $\gamma = -1\ -1\ 0\ 1\ 1$ &\quad $d=6$, $g=4$.
\end{tabular}

Now, for suitable choice of the basis of $\mbox{Pic } X = {\mathbb
Z}^7$, we can
represent these curves by the following divisor classes on $X$:
\begin{description}
\item{a)} $(1;0^6)$ $\gon = 1$, line $L = G_1$ meets $C$ in $2$ points
\item{b)} $(3;1^5,0)$ $\gon = 2$, line $L = F_{16}$ meets $C$ in $2$ points
\item{c)} $(4;2,1^5)$ $\gon = 2$, line $L = G_1$ meets $C$ in $3$ points
\item{d)} $(6;2^6)$ $\gon = 3$, has a trisecant not on $X$.
\end{description}
Thus for each of these curves the gonality is computed by a
multisecant $L$, and in the
first three cases, we can choose $L$ to be a line lying on $X$.  In
th fourth case, we
cannot find a trisecant lying on $X$, so we make one biliaison (i.e.,
replace $C$ by
$C+H$ on $X$, where $H$ is the plane section) and obtain
\begin{description}
\item{d$'$)} $(9;3^6)$.
\end{description}
This last curve $C$ in case d$'$) is a complete intersection of two
cubic surfaces.
This is the exceptional case of Clifford dimension $3$ studied by
Martens \cite{GM1}.  He
shows it has gonality $6$, but Clifford index $3$ given by the linear
system $g_9^3$
giving the embedding in ${\mathbb P}^3$.  Now this curve does have a
trisecant $L$ on
$X$, say $E_1$, and this line computes the gonality.

To prove our result, we use the fact that every smooth ACM curve on
the cubic surface
$X$ is obtained from one of a), b), c), d) by biliaison on the surface
$X$.  These curves
all have gonality computed by multisecants.  We prove our result
then, by induction on
the degree.  Our induction statement is the stronger claim that if
$C$ is sufficiently
general, then there is a multisecant $L$ on $X$ of order $k$, such
that the gonality of
$C$ is $d-k$.  We begin the induction with cases a), b), c), d$'$).

For the induction step, suppose that $C$ on $X$ has degree $d$, a
multisecant $L$ on $X$
of order $k$, and gonality $d-k$.  We take $C'$ a general member of
the linear system
$|C+H|$, and let it specialize to $C_1 = C$ union $C_2 = H$.  Then $s
= C_1 \cdot C_2 =
d$, and $\gon(H) = 2$, since $H$ is a plane cubic curve.  Then by
$(3.1)$, $\gon(C') \ge
\min\{d,d-k+2\}$.  Since $k \ge 2$ in all our starting cases, we conclude
$\gon(C')
\ge d-k+2$.  On the other hand, the degree of $C'$ is $d+3$, and $C'
\cdot L = k+1$.  So
the pencil of planes through $L$ cuts out a $g^1_{d-k+2}$ and we find
$\gon(C') = d-k+2$
as required.

In particular, for the complete intersection curves on $X$, we
recover a weak form of
Basili's result \cite{Bas}, since our proof works only for
sufficiently general curves.

\section{Conclusion}

For curves on quartic surfaces, Farkas \cite{F} has shown that the
gonality of some
special classes of curves is computed by $4$-secants.  His method
does not cover all ACM
curves on quartic surfaces, because he always assumes the surface
contains no rational
and no elliptic curves.

If we apply the methods of this paper to ACM curves lying on surfaces
of degree four and
higher, we obtain only an inequality for the gonality, not an exact
figure, and so we
are unable to answer Question $2.11$ in general.  If the surface
contains a line, and
the curve is either in the biliaison class of the line, or residual
to the line, then
the line becomes a multisecant of high order that computes the
gonality.  This case was
also observed by Paoletti \cite{P}.

To make further progress on Question $2.11$ will require some other technique.

\end{document}